\documentclass[11pt]{article}
\newcommand{\version}{February 14, 2004}

\usepackage{amsthm,amsfonts,amsmath}
\swapnumbers
\theoremstyle{plain}

\newtheorem{thm}{THEOREM}[section]
\newtheorem{lm}[thm]{LEMMA}

\newtheorem{conjecture}[thm]{Conjecture}
\theoremstyle{definition}

\theoremstyle{remark}
\newcommand{\upchi}{\raise1pt\hbox{$\chi$}}
\newcommand{\R}{{\mathord{\mathbb R}}}

\newcommand{\tr}{\mathrm{Tr}}

\renewcommand{\|}{{\Vert}}
\numberwithin{equation}{section}
\begin{document}
                                                                                     
\markboth{\scriptsize{CL \version}}{\scriptsize{CL \version}}

\title{\bf{SOME MATRIX REARRANGEMENT INEQUALITIES}}
\author{\vspace{5pt} Eric Carlen$^1$ and
Elliott H. Lieb$^{2}$ \\
\vspace{5pt}\small{$1.$ School of Mathematics, Georgia Tech,
Atlanta, GA 30332}  \\
\vspace{5pt}\small{$2.$ Departments of Mathematics and Physics, Jadwin
Hall,} \\[-6pt]
\small{Princeton University, P.~O.~Box 708, Princeton, NJ
  08544}\\
 }
\date{\version}
\maketitle
\footnotetext                                                                         
[1]{Work partially
supported by U.S. National Science Foundation
grant DMS 03-00349.    }                                                          
\footnotetext
[2]{Work partially
supported by U.S. National Science Foundation
grant PHY 01-39984.\\
\copyright\, 2003 by the authors. This paper may be reproduced, in its
entirety, for non-commercial purposes.}

\begin{center}{\textit{ Dedicated to Professor Roberto Conti}}
\end{center}
                                                                      
\begin{abstract}

  We investigate a rearrangement inequality for pairs of $n\times n$
  matrices: Let $\|A\|_p$ denote $({\rm Tr}(A^*A)^{p/2})^{1/p}$, the
  $C^p$ trace norm of an $n\times n$ matrix $A$. Consider the quantity
  $\|A+B\|_p^p + \|A-B\|_p^p$.  Under certain positivity conditions,
  we show that this is nonincreasing for a natural ``rearrangement''
  of the matrices $A$ and $B$ when $1 \le p \le 2$. We conjecture that
  this is true in general, without any restrictions on $A$ and $B$.
  Were this the case, it would prove the analog of Hanner's inequality
  for $L^p$ function spaces, and would show that the unit ball in $C^p$
  has the exact same moduli of smoothness and convexity as does the unit
  ball in $L^p$ for all $1 < p < \infty$. At present this is known to
  be the case only for $1 < p \le 4/3$, $p =2$, and $p\ge 4$. Several
  other rearrangement inequalities that are of interest in their own
  right are proved as the lemmas used in proving the main results.

\end{abstract}

\section{Introduction} \label{intro}

We prove certain rearrangement inequalities for matrices, the main
results being Theorems \ref{thm1} and \ref{thm2} below.  These
rearrangement inequalities pertain to the non-commutative
version \cite{BCL, jag} of Hanner's inequalities \cite{ohanner}.

Hanner's inequalities for  any $L^p$ function space state that
\cite[Theorem 2.5]{anal}
\begin{equation}\label{hanner}
\|f+g\|_p^p +\|f-g\|_p^p \geq \left( \|f\|_p +  \|g\|_p\right)^p
+ \big| \|f\|_p - \|g\|_p\big|^p  
\end{equation}
for $1\leq p\leq 2$. The inequality reverses for $2\leq p \leq
\infty$.

Now, specialize to the case of $L^p(\R^n)$ with Lebesgue measure,
and let $f$ and $g$ be non  negative functions on $\R^n$. Further, let $f^*$
and $g^*$ denote their respective spherical symmetric decreasing
rearrangements. (See, e.g., \cite{anal} for definitions.)  The
Chiti--Tartar inequality states \cite[Theorem 3.5]{anal} that for any
$p\ge 1$
\begin{equation}\label{chiti}
\|f-g\|_p \ge \|f^* - g^*\|_p\ , 
\end{equation} It is also well known that
\begin{equation}\label{plus}
\|f+g\|_p \le \|f^* + g^*\|_p\ \ .
\end{equation}

One can extend the notion of spherical symmetric decreasing
rearrangement to complex valued functions simply by putting $f^* =
|f|^*$.  Under this extension, (\ref{chiti}) remains valid, but
(\ref{plus}) does not.

To connect (\ref{chiti}) and (\ref{plus}) to (\ref{hanner}), consider
the sum $\|f+g\|_p^p + \|f-g\|_p^p$. It is natural to ask how this
compares with $\|f^*+g^*\|_p^p + \|f^*-g^*\|_p^p$. For $p=2$, the
answer is clear; both quantities reduce to $2(\|f\|_2^2+\|g\|_2^2)$.

For other values of $p$, we have the following, which, like (\ref{chiti}),
does \textit{not require
non-negativity} of $f$ and $g$:

\begin{lm} \label{almgren}
For all $1\le p \leq 2$, and all complex-valued  functions $f$
and $g$ in $L^p(\R^n)$
\begin{equation} \label{functions}
\|f+g\|_p^p + \|f-g\|_p^p \geq \|f^*+g^*\|_p^p + \|f^*-g^*\|_p^p  \ .
\end{equation}
For $p>2$, the inequality reverses.
\end{lm} 

\noindent{\bf Proof:} For non negative $f$ and $g$, this lemma follows 
from a theorem in \cite[Theorem 2.2]{AL}, which states that for
non-negative $f$ and $g$, $\int J(f(x), g(x))dx \geq \int J(f^*(x),
g^*(x))dx $ if ${\displaystyle \frac{\partial^2}{ \partial x\partial
    y}J(x,y) \leq 0}$. It is easy to check that, with $J(x,y) =
|x+y|^p + |x-y|^p\ $, we have ${\displaystyle\frac{\partial^2}{
    \partial x\partial y}J(x,y) \leq 0}$ for $1\leq p\leq 2$ and $\geq
0$ for $2\leq p <\infty$.  This establishes (\ref{functions}) for
non-negative functions.

To complete the proof for complex functions it clearly suffices to
prove that if $f,\ g$ are replaced by $|f|,\ |g|$ then the left side
of (\ref{functions}) decreases (resp. increases) for $1\leq p \leq 2$
(resp. $2\leq p <\infty$). A pointwise inequality suffices for this.
For any two real numbers $a$ and $b$, and any $-1 \le t \le 1$, define
the functions $c(t)$ by
\begin{equation}\label{cdef}
c(t) = (a^2 + b^2 + 2abt)^{p/2} + 
(a^2 + b^2 - 2abt)^{p/2}\ .
\end{equation}
This function is strictly concave for $p<2$, and strictly convex for
$p>2$. Since $c'(0) =0$ in either case, we have that $t=\pm1 $
minimizes $c$ for $p<2$, while $t=\pm 1$ maximizes $c$ for $p>2$.
Taking $a = |f(x)|$, $b = |g(x)|$ and $t = (\overline{f(x)}g(x) +
f(x)\overline{g(x)})/(2|f(x)||g(x)|)$, we see that for $p<2$
\begin{equation} \label{pointwise}
|f(x)- g(x)|^p+|f(x) + g(x)|^p \ge (|f(x)| + |g(x)|)^p +  \big||f(x)|- 
|g(x)|\big|^p\ ,
\end{equation}
and this provides the extension to the complex case. \qed
\medskip

Thus, we have the following situation: Consider all complex-valued
functions $\phi$ and $\gamma$ such that $|\phi|, \ |\gamma|$ are
equimeasurable with $|f|$ and $|g|$, respectively (and, therefore,
have the same right side of (\ref{hanner}) as the $f$,$g$ pair).
Which choice will minimize the left side of (\ref{hanner}) when $1\leq
p \leq 2$ (and maximize it when $2\leq p <\infty$)?

One answer is the pair $f^*$ and $g^*$. It is not the only answer,
however, since any further equimeasurable rearrangement such that the
set of level sets of $f$ and those of $g$ are the same will minimize
the left side of (\ref{hanner}).

A converse to Lemma \ref{almgren} is the following, which says that
the parallelogram identity for $p=2$ holds as an inequality for $p\ne
2$.
\begin{lm} \label{antialmgren}
For all $1\le p \leq 2$, and all complex-valued  functions $f$
and $g$ in $L^p(\R^n)$
\begin{equation} \label{antifunctions}
\|f+g\|_p^p + \|f-g\|_p^p \leq \lim_{y\to \infty}\left(\|f+\tau_y g\|_p^p + 
\|f-\tau_y g\|_p^p\right) 
= 2 \|f\|_p^p + 2  \|g\|_p^p \ ,
\end{equation}
where $\tau_y$ is translation by $y\in \R^n$.
For $p>2$, the inequality reverses.
\end{lm}
The proof follows from the fact that that the function $c$ in (\ref{cdef})
is maximized at $t=0$ for $p< 2$, 
and minimized there for $p>2$ by the same convexity argument.

Our goal here is to extend Lemmas \ref{almgren} and \ref{antialmgren}
to matrices.  This will first require some discussion and notation.

\section{Definitions and Main Theorems}\label{definitions}

Let $A$ be any $n\times n$ matrix. Then $|A| = \sqrt{A^*A}$ and,
for $1\le p < \infty$, $\|A\|_p=\left(\tr |A|\right)^{1/p}$.
The analogue of Hanner's inequality (\ref{hanner}) is

\begin{equation}\label{hannerm}
\|A+B\|_p^p +\|A-B\|_p^p \geq \left( \|A\|_p +  \|B\|_p\right)^p
+ \left( \|A\|_p - \|B\|_p\right)^p   
\end{equation}
\textit{for $1\leq p\leq 2$. The inequality reverses for $2\leq p \leq
\infty$.}

\medskip
Inequality (\ref{hannerm}) was proved in \cite{BCL} for the following cases:
\begin{enumerate}
\item For all $1\leq p \leq 4/3$ and $4\leq p \leq \infty$, and of
  course $p=2$.
\item For all $1\leq p \leq \infty$ if $A+B$ and $A-B$ are positive semidefinite.
\label{allp}
\end{enumerate}

We conjecture that (\ref{hannerm}) holds for all $A,B$.  (In
\cite{BCL} the condition in item \ref{allp} was incorrectly stated for
the case $2\leq p \leq \infty $; we are grateful to C. King for
pointing out this error.)

Let $\sigma_1 \ge \sigma_2 \ge \dots \ge \sigma_n \geq 0$ be the
singular values of an $n\times n$ matrix $A$, i.e., the eigenvalues of
$|A|$. . Let $\Sigma_\uparrow (A)$ and $\Sigma_\downarrow (A)$ be the
$n\times n$ matrices defined by
\begin{equation}
\Sigma_\uparrow (A) = \left[
\begin{matrix}
\sigma_1 & \phantom{\sigma_1} &   \phantom{\sigma_1} & \phantom{\sigma_1} \\
\phantom{\sigma_1} & \sigma_2 &   \phantom{\sigma_1} & \phantom{\sigma_1} \\
\phantom{\sigma_1} & \phantom{\sigma_1} &  \ddots & \phantom{\sigma_1} \\
\phantom{\sigma_1} & \phantom{\sigma_1} & \phantom{\sigma_1} & \sigma_n 
\end{matrix} \right]
\qquad{\rm and}\qquad
\Sigma_\downarrow (A) = \left[
\begin{matrix}
\sigma_n & \phantom{\sigma_1}  & \phantom{\sigma_1} & \phantom{\sigma_1} \\
\phantom{\sigma_1} & \sigma_{n-1} &  \phantom{\sigma_1} & \phantom{\sigma_1} \\
\phantom{\sigma_1} & \phantom{\sigma_1} & \ddots & \phantom{\sigma_1} \\
\phantom{\sigma_1} & \phantom{\sigma_1} & \phantom{\sigma_1} & \sigma_1 
\end{matrix}
\right]\ .
\end{equation}  
\medskip

We note, for later use in the proof of Theorems \ref{thm1} and \ref{thm2}, that 
if $A$ and $B$ are Hermitean and $A>|B|$ then
\begin{equation}\label{monotone}
\Sigma_\uparrow (A) >  \Sigma_\uparrow (B)\qquad  {\mathrm{and}} 
\qquad \Sigma_\downarrow (A) >  \Sigma_\downarrow (B) \  .
\end{equation}

The rearrangement $\Sigma_\downarrow(A)$ is considered in the book of
Horn and Johnson \cite{hj} in the same notation apart from the arrow.
(They only consider the decreasing ordering of the singular values).
In problem 18 in section 3.5 \cite{hj}, a proof is sketched of the
analog of (\ref{chiti}) for matrices: Namely that
$$\|\Sigma_\downarrow(A) - \Sigma_\downarrow(B)\| \le \|A - B\|$$
for any unitarily invariant norm, and hence for the $C^p$ norms in
particular. However, the methods employed there do not seem to be
useful when the direction of the inequality depends on the particular
unitarily invariant norm under consideration, as in the matrix analogs
of Lemmas \ref{almgren} and \ref{antialmgren}, which are the following
conjectures:

\begin{conjecture}\label{conjone}
  For all $1\le p \leq 2$, and all complex-valued $n\times n $
  matrices $A$ and $B$
\begin{equation} \label{matrices}
\|A+B\|_p^p + \|A-B\|_p^p  \ge  \|\Sigma_\uparrow (A)+
\Sigma_\uparrow (B)\|_p^p + 
\|\Sigma_\uparrow (A)-\Sigma_\uparrow (B)\|_p^p \ .
\end{equation}
For $p>2$, the inequality reverses.
\end{conjecture}

\begin{conjecture}\label{conjtwo}
  For all $1\le p \leq 2$, and all complex-valued $n\times n $
  matrices $A$ and $B$
\begin{equation} \label{matricesa}
\|A+B\|_p^p + \|A-B\|_p^p  \leq   \|\Sigma_\uparrow (A)+
\Sigma_\downarrow (B)\|_p^p + 
\|\Sigma_\uparrow (A)-\Sigma_\downarrow (B)\|_p^p \ .
\end{equation}
For $p>2$, the inequality reverses.
\end{conjecture}

Note that Conjecture \ref{conjone}, if true, implies (\ref{hannerm}) in
full generality because it reduces the matrix case to the commutative
case, Theorem \ref{hanner}, namely to diagonal matrices (which are
just functions on $\{1,\ 2,\ \cdots,\ n\}$). We also note that
(\ref{matrices}) holds with the  reverse inequality for $p$ an
even integer and also without restriction on $A$ and $B$. This is some
evidence for the validity of the conjecture.

We can prove the following cases (where $X\geq Y$ means that $X-Y$ 
is positive-semidefinite).
\begin{thm} \label{thm1}
Conjecture \ref{conjone} is true for $1\leq p \leq 2$ 
if $A$ and  $B$ are self adjoint and 
$A \geq B \geq 0$.
\end{thm}

\begin{thm} \label{thm2}
Conjecture \ref{conjtwo} is true  for $1\leq p \leq 2$ 
if $A$ and $ B$ are self adjoint and  $A \geq |B| $.
\end{thm}

While both theorems contain positivity conditions, at least Theorem
\ref{thm2} does not require $B$ to be positive.  The presence of the
positivity conditions in Theorems \ref{thm1} and \ref{thm2} reflects
an important difference between the cases of matrices and functions.
In the case of functions, the simple pointwise inequality
(\ref{pointwise}) sufficed to reduce matters to the consideration of
positive functions.  In the case of matrices, this is not possible:
There is actually an inequality that goes in the direction
\textit{opposite} to (\ref{pointwise}).

\begin{lm} \label{otherway}
  Let $A$ and $B$ be self adjoint $n\times n$ matrices, and suppose
  that $A \ge |B| \ge 0$. Then for $1\le p \le 2$,
$$\tr\left((A+B)^p + (A- B)^p\right)  \le   
\tr\left((A+|B|)^p + (A- |B|)^p \right)\ .$$

\end{lm} 

\medskip
\noindent{\bf Proof:} 
Let $X$ denote the positive part of $B$, and let $Y$ denote the
negative part so that $B = X - Y$ and $|B| = X+Y$.  Define functions
$f(t)$ and $g(t)$ for $0 \le t\le 1$ by
$$f(t) = \frac{1}{ p}\tr\left((A+tB)^p + (A- tB)\right)^p\qquad{\rm
  and}\qquad g(t) = \frac{1}{ p}\tr\left((A+t|B|)^p + (A-
  t|B|)\right)^p\ .$$
Clearly $f(0) = g(0)$, We claim that for each
$t$ with $0 < t \le 1$, $f'(t) > g'(t)$.  To see carry out the
computation that demonstrates this, define the positive semidefinite
matrices
$$Z_1 = A+t(X+Y)\qquad Z_2 = A+ t(X-Y)\qquad Z_3 = A-t(X-Y)\qquad{\rm and}\qquad
Z_4 = A - t(X+Y)\ .$$ Then
\begin{eqnarray} 
&\phantom{.}&f'(t) - g'(t) =\nonumber\\
&\phantom{.}&\tr\left((X+Y)(Z_1^{p-1} -
Z_4^{p-1})\right) - \tr\left((X-Y)(Z_2^{p-1} - Z_3^{p-1})\right) =
\nonumber\\
&\phantom{.}&\tr\left(X(Z_1^{p-1} +Z_3^{p-1} - Z_2^{p-1} - Z_4^{p-1})\right) +
\tr\left(Y(Z_1^{p-1} + Z_2^{p-1} - Z_3^{p-1} - Z_4^{p-1})\right) =
\nonumber\\
&\phantom{.}&\tr\left(X([Z_1^{p-1} - Z_2^{p-1}] +[Z_3^{p-1} - Z_4^{p-1}])\right) +
\tr\left(Y([Z_1^{p-1} - Z_3^{p-1}] + [Z_2^{p-1} - Z_4^{p-1}])\right) \
.\nonumber\end{eqnarray}

Because $0 \le p-1 \le 1$, the operator monotonicity of $(p-1)$st
powers implies that all of the differences in square brackets are
positive. Hence $f'(t) - g'(t) \ge 0$, and so $f(1) \ge g(1)$.  \qed

\medskip

\section{Proof of Theorems \ref{thm1} and \ref{thm2}}

The proofs of both theorems rely on a rearrangement inequality for
alternating products of two positive $n\times n$ matrices $A$ and $B$.
The fact needed in the proofs is that, for integer $s>0$, the quantity
\begin{equation}\label{product}
\tr(BABA\dots BAB) = \tr(B(B^{1/2}AB^{1/2})^s)
\end{equation}
is nonincreasing if we rearrange $A$ and $B$ oppositely, and
nondecreasing if we rearrange $A$ and $B$ similarly. That is, the
quantity in (\ref{product}) does not increase if we replace $A$ by
$\Sigma_\uparrow(A)$ and $B$ by $\Sigma_\downarrow(B)$, and does not
decrease if replace $A$ by $\Sigma_\uparrow(A)$ and $B$ by
$\Sigma_\uparrow(B)$.  The following theorems assert this, and
somewhat more.

\begin{thm}\label{updownthm1}  For any two positive--semidefinite
  $n\times n$ matrices $A$ and $B$, any numbers $r\ge 0$ and $s\ge 1$,
\begin{equation}\label{updown1}
 \tr\left(B^r (B^{1/2}AB^{1/2})^s\right) \ge 
\tr\left((\Sigma_\uparrow (A))^s(\Sigma_\downarrow (B))^{s+r}\right)\ .
\end{equation}
\end{thm}

\begin{thm}\label{updownthm2}  For any two positive--semidefinite
  $n\times n$ matrices $A$ and $B$, any number $r\ge 0$ and any
  integer $s\ge 1$,
\begin{equation}\label{updown2}
\tr\left((\Sigma_\uparrow (A))^s(\Sigma_\uparrow (B))^{s+r}\right)
\ge  \tr\left(B^r (B^{1/2}AB^{1/2})^s\right) \ .
\end{equation}
\end{thm}

Unlike Theorem \ref{updownthm1}, Theorem \ref{updownthm2} requires
that $s$ be an integer. This condition on $s$ would be unnecessary if
a natural generalization of an inequality of Lieb and Thirring
\cite[Appendix B]{ltm} were established, as we explain in an appendix
to this paper, where further trace inequalities are conjectured and
proved. These are closely related to results and a conjecture in our
earlier paper \cite{cl}.

One tool used in the proof of Theorems \ref{updownthm1} and
\ref{updownthm2} is a ``layer cake representation'' for positive
matrices. Let $C$ be any positive $n\times n$ matrix with spectral
decomposition $C = \sum_{i=1}^s\lambda_i u_i u_i^*$ in which the
eigenvalues $\lambda_j$ are arranged in decreasing order. Let
\begin{equation}\label{layer1}
P_j = \sum_{i=1}^j u_i u_i^*\ .
\end{equation}
Then $P_j$ is the
orthogonal projection onto an eigenspace corresponding to the $j$
largest eigenvalues of $C$, and clearly $P_j \subset P_{j+1}$. With
$P_0 = 0$, we have
$$C = \sum_{j=1}^n \lambda_j(P_j - P_{j-1}) = \lambda_n P_n +
\sum_{j=1}^{n-1}(\lambda_j - \lambda_{j+1})P_j\ .$$ Define  $c_j = \lambda_j -
\lambda_{j+1}$  for $1\le j\le n-1$ and  $c_n =
\lambda_n$. Then
\begin{equation}\label{layer2}
C = \sum_{j=1}^n c_jP_j\ .
\end{equation}  
Note that each $c_j$ is non negative, and $\sum_{j=1}^n c_j =
\lambda_1$. \textit{Therefore, if $\|C\|_\infty =1$, i.e., $\lambda_1
  =1$, then (\ref{layer2}) presents $C$ as a convex combination of
  projections}.

\medskip \noindent{\bf Proof of Theorem \ref{updownthm1}:}  
\textit{(Step One: Reduction to the case $r=0$)}
Observe that,
with  $X = B^{1/2}AB^{1/2}$ and $Y = B^{r/s}$, 
$\tr(B^r(B^{1/2}AB^{1/2})^s) = 
\tr(X^sY^s)$. The  inequality of  \cite[Appendix B]{ltm}  asserts that
$\tr(X^sY^s) \ge \tr((Y^{1/2}XY^{1/2})^s)$. Therefore, 
$$\tr(B^r(AB)^s) \ge \tr((B^{(s+r)/2s}AB^{(s+r)/2s})^s)\ .$$

Given the validity of (\ref{updown1}) in the case $r=0$, we have
\begin{eqnarray}
\tr((B^{(s+r)/2s}AB^{(s+r)/2s})^s) &\ge& 
\tr\left((\Sigma_\downarrow (B)^{(s+r)/2s})^{s}
(\Sigma_\uparrow (A))^s(\Sigma_\downarrow (B)^{(s+r)/2s})^{s}\right)\nonumber\\
&=& \tr\left((\Sigma_\uparrow (A))^s(\Sigma_\downarrow 
(B))^{s+r}\right)\ .\end{eqnarray}

\textit{(Step Two: Proof of (\ref{updown1}) for $r=0$)} This is based
on Epstein's concavity theorem \cite{Ep}, and the layer cake
representation (\ref{layer2}). Without loss of generality, we may
suppose that $\|A\| = \|B\| =1$.

Let $C= A^s$, and note that $\|C\| = 1$ as well. Then 
$$\tr(B^{1/2}A B^{1/2})^s = \tr(B^{1/2}C^{1/s}B^{1/2})^s\ .$$
By
Epstein's theorem, $\tr(B^{1/2}C^{1/s}B^{1/2})^s$ is a concave
function of $C$.  Since $\|C\| =1$, the layer cake representation $C =
\sum_{j=1}^n c_jP_j$ is a convex combination of projections, and hence
\begin{equation}\label{step}
\tr(B^{1/2}C^{1/s}B^{1/2})^s \ge
\sum_{j=1}^nc_j\tr(B^{1/2}P_j^{1/s}B^{1/2})^s\end{equation}
Next, since each $P_j$ is an orthogonal projection, 
$P_j^{1/s} = P_j = P_j^{1/2}$. Also, for any two positive 
semidefinite matrices $X$ and $Y$, $\tr((X^{1/2}YX^{1/2})^s) = 
\tr((Y^{1/2}XY^{1/2})^s)$, since, in the positive definite 
case, $X^{1/2}YX^{1/2}$ and
$Y^{1/2}XY^{1/2}$ are similar matrices.
Hence (\ref{step}) becomes
$$\tr(B^{1/2}C^{1/s}B^{1/2})^s \ge \sum_{j=1}^nc_j\tr(P_jBP_j)^s\ .$$
Now do the same thing for $B$: Let $D = B^s$, and let $D =
\sum_{k=1}^n d_k Q_k$ be the layer--cake representation of $D$. Again,
since the largest eigenvalue of $D$ is $1$, this displays $D$ as a
convex combination of projections.  Again applying Epstein's theorem,
and using the fact that each $Q_k^{1/s} = Q_k$, we deduce that for
each $j$,

\begin{equation} 
\tr(P_jBP_j)^n  = \tr (P_jD^{1/s}P_j)^s
 \ge
\sum_{k=1}^nd_k \tr(P_jQ_kP_j)^s
 = \sum_{k=1}^n d_k
\tr(Q_kP_jQ_k)^s\ .
\end{equation}
Now, $Q_kP_jQ_k$ is a positive contraction. A vector $v$ is an
eigenvector of this contraction with eigenvalue $1$ if and only if $v$
belongs to the images of both $P_j$ and $Q_k$. Let $r_{j,k}$ be the
dimension of the intersection of the images of both $P_j$ and $Q_k$.
This is the geometric multiplicity of $1$ as an eigenvalue of
$Q_kP_jQ_k$. Since all of the other eigenvalues are non negative,
$$\tr (Q_kP_jQ_k)^s \ge r_{j,k}$$
for all $s$. (In fact, it converges to this value as $s$ increases).

For any two subspaces ${\cal V}_1$ and ${\cal V}_2$,
$${\rm dim}({\cal V}_1 + {\cal V}_2) + {\rm dim}({\cal V}_1 
\cap {\cal V}_2) =
{\rm dim}({\cal V}_1) + {\rm dim}({\cal V}_2)\ ,$$
and we have
$$ r_{j,k} \ge \max\{\ \tr(P_j) +  \tr(Q_k) - N\ ,\ 0\ \}\ .$$
There is equality in this inequality if we replace $P_j$ by 
$\Sigma_\uparrow(P_j)$, and 
$Q_j$ by $\Sigma_\downarrow(Q_k)$.
Hence, since these matrices commute,
\begin{eqnarray}
 \tr(AB)^s &\ge& \sum_{j,k=1}^n
c_jd_k\tr\left((\Sigma_\uparrow(P_j)
\Sigma_\downarrow(Q_k))^s\right)
=
\tr\left(\left(\sum_{j=1}^n
c_j\Sigma_\uparrow(P_j)\right) \left(\sum_{k=1}^n d_k
\Sigma_\downarrow(Q_k))\right)\right)\nonumber\\
&=&
\tr\left(\Sigma_\uparrow(A^s)\Sigma_\downarrow(B^s)\right)
= \tr\left(\left(\Sigma_\uparrow (A)\Sigma_\downarrow
(B)\right)^s\right)\ .\end{eqnarray} \qed

\medskip
\noindent{\bf Proof of Theorem \ref{updownthm2}}:
Let $N$ denote the integral value of $s$. Notice that $\tr\left(B^r
  (B^{1/2}AB^{1/2})^N\right) = \tr\left(B^r (AB)^N\right)$.  Expand
$A$, $B$ and $B^r$ in their layer cake representations as above.
Taking the trace, we get a linear combination with positive
coefficients of terms such as ${\displaystyle \prod_{j=1}^{2N+1}P_j } $,
where each $P_j$ is a projection coming from
one of the layer cake expansions.  By cyclicity of the trace, we may
assume that $\tr(P_1) = \min\left\{\ \tr(P_j) \ :\ 1\le j \le 2N+1\ 
\right\}$. Then since ${\displaystyle \prod_{j=2}^{2N+1}P_j}$ is a
contraction, if we compute ${\displaystyle
  \tr\left(\prod_{j=1}^{2N+1}P_j\right)}$ in a basis of eigenvectors
of $P_1$, we certainly find that
$$\left|\tr\left(\prod_{j=1}^{2N+1}P_j\right)\right| \le \tr(P_1) 
= \min\left\{\ \tr(P_j) \ :\
1\le j \le 2N+1\ \right\}\ .$$ There is equality in case the $P_j$
all commute, and are all ``nested'' which is what happens if we
replace $A$ and $B$ by $\Sigma_\uparrow (A)$ and $\Sigma_\uparrow (B)$
respectively, or just as well, by $\Sigma_\downarrow (A)$ and
$\Sigma_\downarrow (B)$ respectively. \qed

\medskip 
\noindent{\bf Proof of  Theorem \ref{thm1}:} For any positive matrix
$C$ and any $p$ with $1 < p < 2$,
\begin{equation}\label{intrep}
C^{p} = k_p\, C\, \int_0^\infty \left(\frac{1}{ t} - \frac{1}{
    t+C}\right)t^{p-1}{\rm d}t =
k_p\ \int_0^\infty \left(\frac{C}{ t^2 } - \frac{1}{
    t} + \frac{1}{t+C} \right)t^{p}{\rm d}t \ ,
\end{equation}
where $k_p>0$ is a normalization constant.  We alternately set
$C=A+B>0$ and $C=A-B>0$ in the integrand of (\ref{intrep}) and take
the trace. Since $\tr(A) =\tr( \Sigma_\uparrow(A) )$ and, by
(\ref{monotone}), $\Sigma_\uparrow(A)\pm \Sigma_\uparrow(B)>0$, it
suffices, for our proof, to show that for each $t>0$
\begin{equation}\label{suffice}
\tr\left((t+A+B)^{-1} + (t+A-B)^{-1}\right) 
\ge
\tr\left( (t+\Sigma_\uparrow(A)+\Sigma_\uparrow(B))^{-1} 
+ (t+\Sigma_\uparrow(A)-\Sigma_\uparrow(B))^{-1}\right)
\end{equation}

Let $H$ denote $A+t$. By (\ref{monotone}) we have that $K := H^{-1/2}B H^{-1/2}$ 
satisfies $0<K< 1$. Therefore, it is legitimate to expand
\begin{equation}
(H\pm B)^{-1} = H^{-1/2} (1\pm K)^{-1}H^{-1/2} =  
H^{-1/2}  \sum_{j=0}^\infty (-1)^j (\pm K)^j   H^{-1/2}\ .
\end{equation}
If these two expressions, $\pm$, are added, the left side
of (\ref{suffice}) becomes
\begin{equation} \label{sum}
2 \sum_{j=0}^\infty  \tr (H^{-1} K^{2j}) = 
2 \sum_{j=0}^\infty \tr H^{-1}(B H^{-1})^{2j} \ .
\end{equation}

An expression similar to this is obtained for the right side of
(\ref{suffice}), except that $B$ is replaced by $\Sigma_{\uparrow}(B)$
and $H^{-1}$ is replaced by $\Sigma_{\downarrow}(H^{-1})$, which arises from the 
fact that $ (t+\Sigma_{\uparrow}(A))^{-1} = \Sigma_{\downarrow}((t+A)^{-1})$.
By Theorem (\ref{updownthm1}), this replacement cannot increase each term in
(\ref{sum}). 
\qed

\medskip 
\noindent{\bf Proof of  Theorem \ref{thm2}:}  By Lemma \ref{otherway}, 
we may replace $B$ by $|B|$,
and since $A \ge |B|$, both $A+|B|$ and $A - |B|$ are non negative.
Hence, the integral representation used in the proof of Theorem
\ref{thm1} may be applied. Instead of rearranging oppositely and
applying Theorem \ref{updownthm1}, rearrange similarly, and apply
Theorem \ref{updownthm2}.  Using one additional but obvious fact --
that $\Sigma_\uparrow(B) = \Sigma_\uparrow(|B|)$ -- the theorem is
proved. \qed

\bigskip

\appendix \section{Appendix: Remarks on Theorem \ref{updownthm2}}

We have made use of the inequality  
\begin{equation}\label{liebth}
\tr(Y^{1/2}XY^{1/2})^s \le \tr(X^sY^s)\ ,
\end{equation}
valid for all positive semidefinite $n \times n$ matrices $X$ and $Y$
and all $s\ge 1$.  This inequality was proved in \cite[Appendix
B]{ltm} using Epstein's theorem \cite{Ep}, which asserts the concavity
for all $s\ge 1$ of the function $f_s$, given by
\begin{equation}\label{eps}
f_s(A) =  \tr((B^{1/2}A^{1/s}B^{1/2})^s)
\end{equation}
on the set of positive semidefinite $n\times n$ matrices, where $B$
is some fixed positive semidefinite $n\times n$ matrix.  In \cite{cl},
we conjectured that for $1/2 < s < 1$, $f_s$ is convex.  Indeed,
$f_{1/2}$ {\it is} convex, since $\tr((B^{1/2}A^{2}B^{1/2})^{1/2}) =
\|AB^{1/2}\|_1$.  We used the concavity of $f_s$ for $s \ge 1$ to
prove a Minkowski type inequality for traces, and showed how this
yielded a proof of the strong subadditivity of the quantum mechanical
entropy.  The conjectured convexity would have done the same thing.
 
 Were $f_s$ convex for $1/2 < s < 1$, the kind of proof in \cite{ltm}
 of (\ref{liebth}) would carry over to a proof that for all positive
 semidefinite $n\times n$ matrices $A$ and $B$, and all such $s$,
\begin{equation}\label{rev}
\tr(A^sB^s) \le \tr((B^{1/2}AB^{1/2})^s)
\end{equation}

Using the convexity of $f_{1/2}$ proved above, we now prove (\ref{rev})
for $s = 1/2$.  Let $C = A^{1/2}$, and introduce $g(C)$ and $h(C)$
where $g(C) = \tr(CB^{1/2})$ and $h(C) =
\tr(B^{1/2}C^2B^{1/2})^{1/2}$.  By what has been said above, $h(C) -
g(C)$ is a convex function of $C$.  Let $C^{(d)}$ denote the part of
$C$ that is diagonal in a basis that diagonalizes $B$, and let
$C^{(o)}$ denote the off diagonal part. As in \cite{ltm}, one sees
that
$$k(t) = h(C^{(d)} + t C^{(o)}) - g(C^{(d)} + t C^{(o)})$$
satisifies
$k(0) = k'(0) =0$, and since it is convex, $k(1) \ge 0$. This proves
(\ref{rev}) for $s =1/2$.  Given the conjectured convexity of $f_s$
for $1/2 < s < 1$, this argument would establish (\ref{rev}) for $1/2
< s < 1$.

Another possible generalization of (\ref{liebth}) is the following:
\begin{equation}\label{liebth2} \tr(Y^r(Y^{1/2}XY^{1/2})^s) \le
  \tr(X^sY^{s+r})
\end{equation}
for all positive semidefinite $\times n$ matrices $X$ and $Y$, all
$r>0$, and all $s\ge 1$. Were this true, one could easily remove the
restriction that $s$ be an integer in Theorem \ref{updownthm2}. This
would be true if it were the case that
$$\tilde f_s(A) = \tr(B^r(B^{1/2}A^{1/s}B^{1/2})^s))$$
defined a
concave function of the positive semidefinite matrix $A$ for $s \ge
1$, and $B$ a given positive semidefinite matrix.


\medskip






\begin{thebibliography}{99}
\bibitem{AL} F.J.~Almgren and E.H.~Lieb, \textit{ Symmetric Decreasing
    Rearrangement is Sometimes Continuous}, Jour. Amer. Math. Soc.
  {\bf 2}, 683-773 (1989).

\bibitem{BCL} K.~Ball, E.~Carlen and E.H.~Lieb \textit{Sharp Uniform
    Convexity and Smoothness Inequalities for Trace Norms}, Invent.
  Math. {\bf 115}, 463-482 (1994).
  
  
  
\bibitem{cl} E.~Carlen and E.H.~Lieb \textit{A Minkowski Type Trace
    Inequality and Strong Subadditivity of Quantum Entropy}, Advances
  in the Mathematical Sciences, AMS Translations, {\bf 189} Series 2,
  (1999) 59-68.  Also in \textit{Inequalities, Selecta of Elliott H.
    Lieb} M. Loss, M.B. Ruskai eds., Springer, 2002.

  
  
\bibitem{Ep} H.~Epstein, \textit{Two Theorems of E. Lieb}, Commun.
  Math. Phys. {\bf 31}, 317-322 (1973).

\bibitem{ohanner} O.~Hanner, \textit{ On the uniform convexity of 
$L^p$ and $\ell^p$}, Ark. Math. {\bf 3}, 239-244 (1956).
  
\bibitem{hj}   R.A.~Horn  and C.R.~Johnson,  \textit{Topics in  matrix
    analysis}, Cambridge  University Press, Cambridge. second edition,
  1991.

\bibitem{anal}
E.H.~Lieb and M.~Loss, \textit{Analysis},
Amer. Math. Soc. second edition, 2001.


\bibitem{ltm} E.H.~Lieb and W.~Thirring, \textit{Inequalities for the
    Moments of the eigenvalues of the Schrodinger Hamiltonian and
    Their Realtion to Sobolev Inequalities}, in \textit{Studies in
    Mathematical Physics}, E. Lieb. B. Simon, A. Wightman
  eds.,Princeton University Press, 269--303,1976.

\bibitem{jag} N.~Tomczak-Jaegermann, \textit{The moduli of smootheness and convexity
and Rademacher averages of trace classes $S_p$ ($1\leq p <\infty$)},
Studia Math. {\bf 50}, 163-182 (1974).
                                                                                                     
                                                                                                     





\end{thebibliography}
\end{document}